\documentclass[12pt]{article}
\usepackage{amssymb}
\newcommand{\be}{\begin{equation}}
\newcommand{\ee}{\end{equation}}
\newcommand{\bea}{\begin{eqnarray}}
\newcommand{\eea}{\end{eqnarray}}
\newcommand{\ba}{\begin{array}}
\newcommand{\ea}{\end{array}}

\newcommand{\bc}{\begin{center}}
\newcommand{\ec}{\end{center}}
\newcommand{\ben}{\begin{enumerate}}
\newcommand{\een}{\end{enumerate}}
\newcommand{\bfi}{\begin{figure}}
\newcommand{\efi}{\end{figure}}

\newcommand{\bq}{\begin{quote}}
\newcommand{\eq}{\end{quote}}
\newcommand{\bqu}{\begin{quotation}}
\newcommand{\equ}{\end{quotation}}
\newenvironment{emphit}{\begin{itemize}}{\end{itemize}}
\newcommand{\bemp}{\begin{emphit}}
\newcommand{\eemp}{\end{emphit}}

\newcommand{\bt}{\begin{tabular}}
\newcommand{\et}{\end{tabular}}

\newtheorem{myth}{Theorem}[section]
\newtheorem{mylem}{Lemma}[section]
\newtheorem{mycor}{Corollary}[section]

\newtheorem{mydef}{Definition}[section]

\newtheorem{myrem}{Remark}[section]

\begin{document}
\date{}
\title{On Some Autotopisms Of Non-Steiner Central Loops
}
\author{T. G. Jaiy\'e\d ol\'a\thanks{All correspondence to be addressed to this author.}\\
Department of Mathematics,\\
Obafemi Awolowo University,\\
Ile Ife 220005, Nigeria.\\
jaiyeolatemitope@yahoo.com\\tjayeola@oauife.edu.ng \and
J. O. Ad\'en\'iran \\
Department of Mathematics,\\
University of Abeokuta, \\
Abeokuta 110101, Nigeria.\\
ekenedilichineke@yahoo.com\\
adeniranoj@unaab.edu.ng}\maketitle
\begin{abstract}
An algebraic process for the construction of an autotopism for a
non-Steiner C-loop is described and this is demonstrated with an
example using a known finite C-loop. In every C-loop, two of its
parastrophes are equivalent(equal) to it, if and only if both the
first and second components of the constructed autotopism and its
inverse autotopism are equal to the identity map. Hence, the other
three parastrophes are equivalent(equal) to the C-loop. It is proved
that the set of autotopisms that prevent a C-loop from being a
Steiner loop forms a Steiner triple system.
\end{abstract}

\section{Introduction}
LC-loops, RC-loops and C-loops are loops that satisfy the identities
\begin{displaymath}
(xx)(yz)=(x(xy))z,~(zy)(xx)=z((yx)x)~\textrm{and}~x(y(yz))=((xy)y)z~\textrm{respectively}.
\end{displaymath}
These three types of loops are collectively called central loops. In
the theory of loops, central loops are some of the least studied
loops. They have been studied by Phillips and Vojt\v echovsk\'y
\cite{phi}, \cite{phi1}, \cite{phi2}, Kinyon et. al. \cite{phi3},
\cite{phi4}, \cite{phi5}, Ramamurthi and Solarin \cite{ram}, Fenyves
\cite{fen} and Beg \cite{phd169}, \cite{phd170}. The difficulty in
studying them is as a result of the nature of the identities
defining them when compared with other Bol-Moufang identities. It
can be noticed that in the aforementioned LC identity, the two $x$
variables are consecutively positioned and neither $y$ nor $z$ is
between them. A similarly observation is true in the other two
identities(i.e the RC and C identities). But this observation is not
true in the identities defining Bol loops, Moufang loops and extra
loops. Fenyves \cite{fen} gave three equivalent identities that
define LC-loops, three equivalent identities that define RC-loops
and only one identity that defines C-loops. But recently, Phillips
and Vojt\v echovsk\'y \cite{phi1}, \cite{phi2} gave four equivalent
identities that define LC-loops and four equivalent identities that
define RC-loops. Three of the four identities given by Phillips and
Vojt\v echovsk\'y are the same as the three already given by
Fenyves.

Their basic properties are found in \cite{phi}, \cite{ram},
\cite{fen} and \cite{den}. The left and right translation maps on
the loop $(L,\cdot )$ denoted by $L_x:L\to L$ and $R_x:L\to L$ and
defined as $yL_x=xy$ and $yR_x=yx$ respectively $\forall~x,y\in L$
are bijections. $L$ is said to be left alternative and right
alternative if
\begin{displaymath}
x\cdot xy=x^2y~\textrm{and}~yx\cdot
x=yx^2~\textrm{respectively}~\forall~x,y\in L
\end{displaymath}
Thus, $L$ is said to be alternative if it is both left and right
alternative. $L$ is called a Steiner loop if and only if
\begin{displaymath}
x^2=e~, ~yx\cdot x=y~\textrm{and}~xy=yx~\forall ~x, y\in L.
\end{displaymath}
The set $S(L, \cdot )$ of all bijections in a loop $(L,\cdot )$
forms a group called the permutation group of the loop $(L,\cdot )$.
The triple $(U, V, W)$ such that $U, V, W\in S(L, \cdot )$ is called
an autotopism of $L$ if and only if
\begin{displaymath}
xU\cdot yV=(x\cdot y)W~\forall ~x, y\in L.
\end{displaymath}
The group of autotopisms of $L$ is denoted by $Aut(L, \cdot )$.

An algebraic process for the construction of an autotopism for a
non-Steiner C-loop is described and this is demonstrated with an
example using a known finite C-loop. In every C-loop, two of its
parastrophes are equivalent(equal) to it, if and only if both the
first and second components of the constructed autotopism and its
inverse autotopism are equal to the identity map. Hence, the other
three parastrophes are equivalent(equal) to the C-loop. It is proved
that the set of autotopisms that prevent a C-loop from being a
Steiner loop forms a Steiner triple system.

\begin{mydef}\label{par:loop}(\cite{den}, Page~65)

Let $(L,\theta )$ be a quasigroup. The 5 parastrophes or conjugates or adjugates of $(L,\theta )$ are quasigroups
whose binary operations $\theta^*~,~\theta^{-1}~,~{}^{-1}\theta~,~(\theta^{-1})^*~,~({}^{-1}\theta )^*$
defined on $L$ are given by :
\begin{description}
\item[(a)]
\begin{displaymath}
(L,\theta^*)~:~y\theta ^*x=z\Leftrightarrow x\theta
y=z~\forall~x,y,z\in L.
\end{displaymath}
\item[(b)]
\begin{displaymath}
(L,\theta^{-1})~:~x\theta ^{-1}z=y\Leftrightarrow x\theta
y=z~\forall~x,y,z\in L.
\end{displaymath}
\item[(c)]
\begin{displaymath} (L,{}^{-1}\theta )~:~ z~{}^{-1}\theta
y=x\Leftrightarrow x\theta y=z~\forall~x,y,z\in L.
\end{displaymath}
\item[(d)]
\begin{displaymath}
\bigg(L,\big(\theta ^{-1}\big)^*\bigg)~:~z\big(\theta
^{-1}\big)^*x=y\Leftrightarrow x\theta y=z~\forall~x,y,z\in L.
\end{displaymath}
\item[(e)]
\begin{displaymath}
\bigg(L,\big({}^{-1}\theta \big)^*\bigg)~:~y\big({}^{-1}\theta
\big)^*z=x\Leftrightarrow x\theta y=z~\forall~x,y,z\in L.
\end{displaymath}
\end{description}
\end{mydef}

\begin{myrem}
As it can be seen in Definition~\ref{par:loop}, every quasigroup
$(L,\cdot )$ belongs to a set of 6 quasigroups, called adjugates by
Fisher and Yates \cite{phd72}, conjugates by Stein \cite{phd73},
\cite{phd76} and Belousov \cite{phd77} and parastrophes by Sade
\cite{phd74}. They have been studied by Artzy \cite{phd71}, Lindner
and Steedley \cite{phd171} and a detailed study on them can be found
in \cite{pfl}, \cite{phd39} and \cite{den}. The most recent studies
of the parastrophes of a quasigroup(loop) are by Duplak
\cite{phd154}, Shchukin and Gushan \cite{phd172}, Frank, Bennett and
Zhang \cite{phd173}.
\end{myrem}

\begin{mydef}\label{st:sys}
A Steiner triple system~(S.T.S. for short) $(Q,3)$ on a set $Q$ is a set of
unordered triples $\{a,b,c\}\in (Q,3)$ such that
\begin{description}
\item[(i)] $a,b,c$ are distinct elements of $Q$,
\item[(ii)] to any $a,b\in Q$ such that $a\ne b$ there exists a unique triple
$\{a,b,c\}\in (Q,3)$.
\end{description}
\end{mydef}

\begin{myrem}\label{sts:order}
It is proved in \cite{den} and stated in \cite{pfl} that if $\big|(Q,3)\big|=r$
where $(Q,3)$ is as defined in Definition~\ref{st:sys}, then, $r\equiv 1\bmod 6$ or $r\equiv 3\bmod 6$.
\end{myrem}

\begin{mydef}\label{st:sys1}
Let $(G,\ast )$ and $(H,\star )$ be two distinct groupoids. $(G,\ast
)$ and $(H,\star )$ are said to be equivalent or equal, written as
$(G,\ast )\equiv (H,\star )$ or $(G,\ast )=(H,\star )$ respectively,
if $G=H$ and $'\ast' ='\star'$. That is,  $(G,\ast )$ and $(H,\star
)$ are the same.
\end{mydef}

\section{Autotopisms of Central Loops}

\begin{myth}\label{lc:auto}
A loop $L$ is an LC-loop $\Leftrightarrow (L_x^2,I,L_x^2)\in Aut(L)~\forall~x\in L$.
\end{myth}
{\bf Proof}\\ Let $L$ be an LC-loop $\Leftrightarrow (x\cdot xy)z=(xx)(yz)\Leftrightarrow(x\cdot xy)z=x(x\cdot yz)$ by \cite{den}
$\Leftrightarrow(L_x^2,I,L_x^2)\in Aut(L)~\forall~x\in L$.

\begin{myth}\label{rc:auto}
A loop $L$ is an RC-loop $\Leftrightarrow(I,R_x^2,R_x^2)\in Aut(L)~\forall~x\in L$.
\end{myth}
{\bf Proof}\\ Let $L$ be an RC-loop, then $z(yx\cdot x)=zy\cdot xx~\Leftrightarrow y(yx\cdot x)=(zy\cdot x)x$ by \cite{den}
$\Leftrightarrow(I,R_x^2,R_x^2)\in Aut(L)~\forall~x\in L$.

\begin{mylem}\label{c:auto2}
Let $L$ be a C-loop. Then for each $(A,B,C)\in Aut(L,\cdot)$, there
exists a unique pair of $(S_1,T_1,\mathcal{R}_1),
(S_2,T_2,\mathcal{R}_2) \in Aut(L,\cdot)$ for each $x\in L$ such
that $L_x^2=S_2^{-1}S_1, R_x^2=T_1^{-1}T_2,
R_x^{-2}L_x^2=\mathcal{R}_2^{-1}\mathcal{R}_1,
\mathcal{R}_1^{-1}\mathcal{R}_2T_2^{-1}T_1S_2^{-1}S_1=I$.
\end{mylem}
{\bf Proof}\\ By Theorem~\ref{lc:auto} and Theorem~\ref{rc:auto};
\begin{displaymath}
(S_1,T_1,\mathcal{R}_1)=(A,B,C)(L_x^2,I,L_x^2)\in Aut(L)
\end{displaymath}
\begin{displaymath}
(S_2,T_2,\mathcal{R}_2)=(A,B,C)(I,R_x^2,R_x^2)\in Aut(L).
\end{displaymath}
Hence, the conditions hold although the identities do not depend on
$(A,B,C)$, but the uniqueness does.

\begin{myth}\label{c:auto2comp}
Let $L$ be a C-loop and let there exist a unique pair of autotopisms
$(S_1,T_1,\mathcal{R}_1), (S_2,T_2,\mathcal{R}_2)$ such that the
conditions $L_x^2=S_2^{-1}S_1, R_x^2=T_1^{-1}T_2$ and
$R_x^{-2}L_x^2=\mathcal{R}_2^{-1}\mathcal{R}_1$ hold for each fixed
$x\in L$. If $\alpha_1=S_1^{-1}, \alpha_2=S_2^{-1},
\beta_1=T_1^{-1},\beta_2=T_2^{-1}, \gamma_1=\mathcal{R}_1^{-1}$ and
$\gamma_2=\mathcal{R}_2^{-1}$, then:
\begin{displaymath}
(x^2y)\alpha_1 =y\alpha_2~~,~~e\alpha_1
=x^{-2}\alpha_2~~,~~x^{m+2}\alpha_1 =x^m\alpha_2
\end{displaymath}
\begin{displaymath}
(yx^2)\beta_2 =y\beta_1~~,~~e\beta_2
=x^{-2}\beta_1~~,~~x^{m+2}\beta_2 =x^m\beta_1
\end{displaymath}
\begin{displaymath}
(x^2yx^{-2})\gamma_1 =y\gamma_2~~,~~e\gamma_1
=e\gamma_2~~,~~x^m\gamma_1 =x^m\gamma_2
\end{displaymath}
for all $m\in \mathbb{Z}$ and $x,y\in L$.
\end{myth}
{\bf Proof}\\
From Lemma~\ref{c:auto2}:
\begin{displaymath}
L_x^2=S_2^{-1}S_1,R_x^2=T_1^{-1}T_2,R_x^{-2}L_x^2=\mathcal{R}_2^{-1}\mathcal{R}_1.
\end{displaymath}
Keeping in mind that a C-loop is power associative and nuclear
square, we have the following proofs.
\begin{enumerate}
\item $L_x^2=S_2^{-1}S_1\Rightarrow yL_x^2=yS_2^{-1}S_1~\forall~y\in L\Rightarrow
yL_{x^2}=yS_2^{-1}S_1\Rightarrow x^2y=yS_2^{-1}S_1\Rightarrow (x^2y)S_1^{-1}=yS_2^{-1}
\Rightarrow x^2y\alpha_1=y\alpha_2$.

Let $y=x^{-2}~ ;~ x^2x^{-2}\alpha_1=x^{-2}\alpha_2\Rightarrow e\alpha_1=x^{-2}\alpha_2$.

Let $y=x^m~ ;~ x^2y\alpha_1=x^2x^m\alpha_1=x^m\alpha_2\Rightarrow x^{m+2}\alpha_1=x^m\alpha_2$.
\item $R_x^2=T_1^{-1}T_2\Rightarrow yR_x^2=yT_1^{-1}T_2~\forall~y\in L\Rightarrow
yx^2=yT_1^{-1}T_2\Rightarrow yx^2T_2^{-1}=yT_1^{-1}\Rightarrow
yx^2\beta=y\beta_1$.

Let $y=x^{-2}~ ;~ yx^2\beta_2=x^{-2}x^2\beta_2=e\beta_2=x^{-2}\beta_1\Rightarrow e\beta_2=x^{-2}\beta_1$.

Let $y=x^m~ ;~ yx^2\beta_2=x^mx^2\beta_2=x^{m+2}\beta_2=x^m\beta_1\Rightarrow x^{m+2}\beta_2=x^m\beta_1$.
\item $R_x^{-2}L_x^2=\mathcal{R}_2^{-1}\mathcal{R}_1\Rightarrow yR_x^{-2}L_x^2=y\mathcal{R}_2^{-1}\mathcal{R}_1~\forall~y\in L\Rightarrow
x^2yx^{-2}=y\mathcal{R}_2^{-1}\mathcal{R}_1\Rightarrow
(x^2yx^{-2})\mathcal{R}_1^{-1}=y\mathcal{R}_2^{-1}\Rightarrow
(x^2yx^{-2})\gamma_1=y\gamma_2$.

Let $y=e~ ;~ (x^2yx^{-2})\gamma_1=(x^2ex^{-2})\gamma_1=(x^2x^{-2})\gamma_1=e\gamma_1
=e\gamma_2\Rightarrow e\gamma_1=e\gamma_2$.

Let $y=x^m~ ;~ (x^2yx^{-2})\gamma_1=(x^2x^mx^{-2})\gamma_1=x^{2+m-2}\gamma_1=x^m\gamma_1
\Rightarrow x^m\gamma_2\Rightarrow x^m\gamma_1=x^m\gamma_2$.
\end{enumerate}

\begin{mycor}\label{c:autotopism}
Let $L$ be a C-loop. An autotopism of $L$ can be constructed if there exists at least an $x\in L$ such that $x^2\neq e$.
The inverse can also be constructed.
\end{mycor}
{\bf Proof}\\ We need Lemma~\ref{c:auto2} and
Theorem~\ref{c:auto2comp}. If $x^2=e$, then the autotopism is
trivial. Since $L$ is a C-loop, using Lemma~\ref{c:auto2} and
Theorem~\ref{c:auto2comp}, it will be noticed that $(\alpha_1
S_2,\beta_1 T_2,\gamma_1 \mathcal{R}_2)\in Aut(L)$ and $(\alpha_2
S_1,\beta_2 T_1,\gamma_2 \mathcal{R}_1)=(\alpha_1 S_2,\beta_1
T_2,\gamma_1 \mathcal{R}_2)^{-1}$. Hence the proof.

\begin{myrem}If $x, y\in L, x\neq y$ such that $x^2=y^2\neq e$, then $x$ and $y$ will generate the same autotopism.
If $x, y\in L$ such that $x^2y^2=e$, then the autotopism generated by $x$ is the inverse of that generated by $y$.
\end{myrem}

\section{C-loops and Steiner loops}
In \cite{phi}, it was shown that every Steiner loop is a C-loop and
Steiner loops are exactly inverse property loops of exponent two.
Hence generally, C-loops are not Steiner loops. Recall that Steiner
loops are totally symmetric loops, whence all parastrophes are
equivalent to them. In this section, for a loop $(L,\cdot )$, if the
triple $(U,V,W)\in Aut(L,\cdot )$ then, $U$, $V$ and $W$ will be
referred to as the first, second and third components of the
autotopism $(U,V,W)$. In the last section, the autotopisms
$(S_1,T_1,\mathcal{R}_1)$ and $(S_2,T_2,\mathcal{R}_2)$ which shall
be referred to as CS-autotopisms were used to construct the
autotopisms $(\alpha_1 S_2,\beta_1 T_2,\gamma_1 \mathcal{R}_2)$ and
$(\alpha_2 S_1,\beta_2 T_1,\gamma_2 \mathcal{R}_1)$. These four
autotopisms are useful to us in this section. Particularly, the
first component $\alpha_1 S_2$ and the second component $\beta_2
T_1$ are of paramount interest.

\begin{myth}
In every C-loop $L$, two of the parastrophes of $L$ are
equivalent(equal) to $L$, if  and only if both the first and second
components of the autotopisms
\begin{displaymath}
(\alpha_1 S_2,\beta_1 T_2,\gamma_1
\mathcal{R}_2)~\textrm{and}~(\alpha_2 S_1,\beta_2 T_1,\gamma_2
\mathcal{R}_1)~\textrm{respectively}
\end{displaymath}
are equal to the identity map. Hence, the other three parastrophes
are equivalent(equal) to $L$.
\end{myth}
{\bf Proof}\\
Using Theorem~\ref{c:auto2comp}; $(x^2y)\alpha_1=y\alpha_2$ and $(yx^2)\beta_2=y\beta_1\Rightarrow
(xx\cdot y)\alpha_1=y\alpha_2$ and $(y\cdot xx)\beta_2=y\beta_1\Rightarrow
(x\cdot xy)\alpha_1=y\alpha_2$ and $(yx\cdot x)\beta_2=y\beta_1\Rightarrow
(x\cdot xy)\alpha_1\alpha_2^{-1}=y$ and $(yx\cdot x)\beta_2\beta_1^{-1}=y\Rightarrow
(x\cdot xy)\alpha_1S_2=y$ and $(yx\cdot x)\beta_2T_1=y$.

Let $x\cdot y=z$, then $(x\cdot z)\alpha_1S_2=y$. Let $x\theta
y=z$($\theta$ replaces $\cdot$). Thus, $\alpha_1S_2 : L\times L\to
L$ is defined by $(x,z)\alpha_1S_2=(x\cdot
z)\alpha_1S_2=y=x\theta^{-1}z\Leftrightarrow x\theta y=z$.
$(L,\theta^{-1} )$ is a parastrophe of $(L,\theta )=(L,\cdot )$ by
Definition~\ref{par:loop}.
If $\alpha_1S_2=I$, then $(L,\theta )\equiv(L,\theta^{-1} )$ i.e $(L,\theta )=(L,\theta^{-1} )$.\\

Let $y\cdot x=t$, then $(t\cdot x)\beta_2T_1=y$. Thus, $\beta_2T_1 :
L\times L\to L~$ is defined by $(t,x)\beta_2T_1=(t\cdot
x)\beta_2T_1=y=t~{}^{-1}\theta~x\Leftrightarrow y\theta
x=t$($\theta$ replaces $\cdot$). $(L,{}^{-1}\theta )$ is a
parastrophe of $(L,\theta )=(L,\cdot )$ by
Definition~\ref{par:loop}. If $\beta_2T_1=I$, then $(L,\theta
)\equiv(L,{}^{-1}\theta )$
i.e $(L,\theta )=(L,{}^{-1}\theta )$.\\

Conversely, assume that $(L,\theta )\equiv(L,\theta^{-1} )$ and
$(L,\theta )\equiv(L,{}^{-1}\theta )$ i.e $(L,\theta
)=(L,\theta^{-1} )$ and $(L,\theta )=(L,{}^{-1}\theta )$ where
$(L,\theta^{-1} )$ and $(L,{}^{-1}\theta)$ are as defined in
Definition~\ref{par:loop}. Recall that ; $(x\cdot xy)\alpha_1S_2=y$
and $(yx\cdot x)\beta_2T_1=y$. Hence, if $z=x\cdot y$ and $t=y\cdot
x$ then ; $(x\cdot xy)\alpha_1S_2=y$ and $(yx\cdot
x)\beta_2T_1=y\Rightarrow (x\cdot z)\alpha_1S_2=y$ and $(t\cdot
x)\beta_2T_1=y\Rightarrow (x\theta z)\alpha_1S_2=y$ and $(t\theta
x)\beta_2T_1=y\Rightarrow (x\theta z)\alpha_1S_2=x\theta^{-1} z$ and
$(t\theta x)\beta_2T_1=t{}^{-1}\theta x\Rightarrow
\alpha_1S_2=I,\beta_2T_1=I$ because $(L,\theta )\equiv(L,\theta^{-1}
)$ and $(L,\theta )\equiv(L,{}^{-1}\theta )$.

The proof of the last part is as follows.
Consider the definitions of the other three parastrophes in Definition~\ref{par:loop}.
\begin{displaymath}
\bigg(L,({}^{-1}\theta)^*\bigg)=\bigg(L,{}^{-1}\Big(\theta^{-1}\Big)
\bigg)
\end{displaymath}
\begin{displaymath}
=\Big\{x,y,z\in L~ :~ y~
{}^{-1}\Big(\theta^{-1}\Big)~z=x~\Leftrightarrow ~x\theta y=z~
\textrm{where}~ \theta~ \textrm{replaces}~ \cdot\Big\},
\end{displaymath}
\begin{displaymath}
\bigg(L,(\theta
^{-1})^*\bigg)=\bigg(L,\Big({}^{-1}\theta\Big)^{-1}\bigg)
\end{displaymath}
\begin{displaymath}
=\Big\{x,y,z\in L~ :~
z\Big({}^{-1}\theta\Big)^{-1}~x=y~\Leftrightarrow ~x\theta y=z~
\textrm{where}~ \theta~ \textrm{replaces}~\cdot\Big\},
\end{displaymath}
\begin{displaymath}
\bigg(L,\theta
^*\bigg)=\bigg(L,\Big({}^{-1}\big(\theta^{-1}\big)\Big)^{-1} \bigg)
\end{displaymath}
\begin{displaymath}
=\Big\{x,y,z\in L~ :~ y~
{}^{-1}\Big(\theta^{-1}\Big)~x=z~\Leftrightarrow ~x\theta y=z~
\textrm{where}~ \theta~ \textrm{replaces}~ \cdot\Big\}.
\end{displaymath}
From the definitions above, it can be seen that the other three
parastrophes can be derived from the first two. Hence, these three
are equivalent to $L$ since the first two are equivalent to $L$ by
the first part.

\begin{mycor}Both the first and second components of the autotopisms
\begin{displaymath}
(\alpha_1 S_2,\beta_1 T_2,\gamma_1
\mathcal{R}_2)~\textrm{and}~(\alpha_2 S_1,\beta_2 T_1,\gamma_2
\mathcal{R}_1)~\textrm{respectively}
\end{displaymath}
of a C-loop $(L,\cdot )$ are equal to the identity map if and only
if $L$ is a Steiner loop.
\end{mycor}
{\bf Proof}\\ Let $\alpha_1S_2=I, \beta_2T_1=I$ then $S_1=S_2,
T_1=T_2$. Whence, $x^2=e$ and $x(xy)=y, (yx)x=y$ since a C-loop is
alternative. This proves that $L$ is a Steiner loop.

Conversely, if $L$ is a Steiner loop, then $L$ is of exponent 2.
Recall that $\alpha_1S_2=L_x^{-2}$ and $\beta_2T_1=R_x^{-2}$.
Hence, $\alpha_1S_2=L_e$ and $\beta_2T_1=R_e\Rightarrow
\alpha_1S_2=I$ and $\beta_2T_1=I$.

\begin{myrem}
The result above generalizes the fact that Steiner loops are exactly
C-loops that are of exponent 2.
\end{myrem}

\begin{myth}Let
\begin{displaymath}
Q=\Bigg\{(S_i, T_i, \mathcal{R}_i), (S_{i+1}, T_{i+1},
\mathcal{R}_{i+1}), (S_iS_{i+1}, T_iT_{i+1},
\mathcal{R}_i\mathcal{R}_{i+1})\in Aut(L)~\Big\vert ~
\end{displaymath}
\begin{displaymath}
(S_i, T_i, \mathcal{R}_i)=(A, B, C)(L_x^2, I, L_x^2), (S_{i+1},
T_{i+1}, \mathcal{R}_{i+1})=(A, B, C)(I, R_x^2, R_x^2),
\end{displaymath}
\begin{displaymath}
(A, B, C)\in Aut(L, \cdot)\Bigg\}_{i\in\mathbb{N}}
\end{displaymath}
be a set of CS-autotopisms of a non-Steiner C-loop. Define
\begin{displaymath}
(Q, 3)=\Bigg\{\bigg\{(S_i, T_i, \mathcal{R}_i), (S_{i+1}, T_{i+1},
\mathcal{R}_{i+1}), (S_iS_{i+1}, T_iT_{i+1},
\mathcal{R}_i\mathcal{R}_{i+1})\bigg\}~\Bigg\vert~i\in
\mathbb{N}\Bigg\}.
\end{displaymath}
Then, (Q, 3) is a Steiner triple system.
\end{myth}
{\bf Proof}\\
We shall show that Definition~\ref{st:sys} is true for $(Q, 3)$.
\begin{enumerate}
\item
\begin{displaymath}
\big(S_i, T_i, \mathcal{R}_i\big), \big(S_{i+1}, T_{i+1},
\mathcal{R}_{i+1}\big), \big(S_iS_{i+1}, T_iT_{i+1},
\mathcal{R}_i\mathcal{R}_{i+1}\big)
\end{displaymath}
are distinct elements of $Q~\forall~i\in \mathbb{N}.$
\item For any $\big(S_i, T_i, \mathcal{R}_i\big), \big(S_{i+1}, T_{i+1}, \mathcal{R}_{i+1}\big)\in Q$, $\big(S_i, T_i, \mathcal{R}_i\big)\not=\big(S_{i+1}, T_{i+1}, \mathcal{R}_{i+1}\big)$ or else $L$ will become a Steiner loop.
There exists a unique autotopism $\big(S_iS_{i+1}, T_iT_{i+1},
\mathcal{R}_i\mathcal{R}_{i+1}\big)\in AUT(L)~\ni$
\begin{displaymath}
\{\big(S_i, T_i, \mathcal{R}_i\big), \big(S_{i+1}, T_{i+1},
\mathcal{R}_{i+1}\big), \big(S_iS_{i+1}, T_iT_{i+1},
\mathcal{R}_i\mathcal{R}_{i+1}\big)\}\in (Q,3)
\end{displaymath}
is distinct.
\end{enumerate}

Thus, by Definition~\ref{st:sys}, $(Q, 3)$ is a Steiner triple
system.

\subsection{Construction}
Let  us now consider the C-loop of order 12 whose multitplication
table is shown in Table~\ref{table:c}.
\begin{table}[!hbp]
\begin{center}
\begin{tabular}{|c||c|c|c|c|c|c|c|c|c|c|c|c|}
\hline
$\cdot $ & 0 & 1 & 2 & 3 & 4 & 5 & 6 & 7 & 8 & 9 & 10 & 11 \\
\hline  \hline
0 & 0 & 1 & 2 & 3 & 4 & 5 & 6 & 7 & 8 & 9 & 10 & 11 \\
\hline
1 & 1 & 2 & 0 & 4 & 5 & 3 & 7 & 8 & 6 & 10 & 11 & 9 \\
\hline
2 & 2 & 0 & 1 & 5 & 3 & 4 & 8 & 6 & 7 & 11 & 9 & 10 \\
\hline
3 & 3 & 4 & 5 & 0 & 1 & 2 & 9 & 10 & 11 & 6 & 7 & 8 \\
\hline
4 & 4 & 5 & 3 & 1 & 2 & 0 & 10 & 11 & 9 & 7 & 8 & 6 \\
\hline
5 & 5 & 3 & 4 & 2 & 0 & 1 & 11 & 9 & 10 & 8 & 6 & 7 \\
\hline
6 & 6 & 7 & 8 & 10 & 11 & 9 & 0 & 1 & 2 & 5 & 3 & 4 \\
\hline
7 & 7 & 8 & 6 & 11 & 9 & 10 & 1 & 2 & 0 & 3 & 4 & 5 \\
\hline
8 & 8 & 6 & 7 & 9 & 10 & 11 & 2 & 0 & 1 & 4 & 5 & 3 \\
\hline
9 & 9 & 10 & 11 & 8 & 6 & 7 & 3 & 4 & 5 & 2 & 0 & 1 \\
\hline
10 & 10 & 11 & 9 & 6 & 7 & 8 & 4 & 5 & 3 & 0 & 1 & 2 \\
\hline
11 & 11 & 9 & 10 & 7 & 8 & 6 & 5 & 3 & 4 & 1 & 2 & 0 \\
\hline
\end{tabular}
\end{center}
\caption{A non-associative C-loop of order 12}\label{table:c}
\end{table}

The construction of an autotopism of the finite C-loop whose
bordered multiplication table is shown in Table~\ref{table:c} is now
given below. $\Pi_\rho$ denotes the right representation set of the
loop.

By Theorem~\ref{c:auto2comp};
\begin{displaymath}
(x^2y)\alpha_1S_2 =y,
\end{displaymath}
\begin{displaymath}
yx^2 =y\beta_1T_2,
\end{displaymath}
\begin{displaymath}
(x^2yx^{-2})\gamma_1\mathcal{R}_2 =y.
\end{displaymath}
Consider
\begin{displaymath}
(x^2y)\alpha_1S_2 =y
\end{displaymath}
and fix $x=4$.
\begin{displaymath}
\textrm{Let}\qquad y=0, \Big(4^2\cdot 0\Big)\alpha_1S_2=0\Rightarrow 2\alpha_1S_2=0.
\end{displaymath}
\begin{displaymath}
\textrm{Let}\qquad y=1, \Big(4^2\cdot 1\Big)\alpha_1S_2=1\Rightarrow 0\alpha_1S_2=1.
\end{displaymath}
\begin{displaymath}
\textrm{Let}\qquad y=2, \Big(4^2\cdot 2\Big)\alpha_1S_2=2\Rightarrow 1\alpha_1S_2=2.
\end{displaymath}
\begin{displaymath}
\textrm{Let}\qquad y=3, \Big(4^2\cdot 3\Big)\alpha_1S_2=3\Rightarrow 5\alpha_1S_2=3.
\end{displaymath}
\begin{displaymath}
\textrm{Let}\qquad y=4, \Big(4^2\cdot 4\Big)\alpha_1S_2=4\Rightarrow
3\alpha_1S_2=4.
\end{displaymath}
\begin{displaymath}
\textrm{Let}\qquad y=5, \Big(4^2\cdot 5\Big)\alpha_1S_2=5\Rightarrow 4\alpha_1S_2=5.
\end{displaymath}
\begin{displaymath}
\textrm{Let}\qquad y=6, \Big(4^2\cdot 6\Big)\alpha_1S_2=6\Rightarrow 8\alpha_1S_2=6.
\end{displaymath}
\begin{displaymath}
\textrm{Let}\qquad y=7, \Big(4^2\cdot 7\Big)\alpha_1S_2=7\Rightarrow 6\alpha_1S_2=7.
\end{displaymath}
\begin{displaymath}
\textrm{Let}\qquad y=8, \Big(4^2\cdot 8\Big)\alpha_1S_2=8\Rightarrow 7\alpha_1S_2=8.
\end{displaymath}
\begin{displaymath}
\textrm{Let}\qquad y=9, \Big(4^2\cdot 9\Big)\alpha_1S_2=9\Rightarrow 11\alpha_1S_2=9.
\end{displaymath}
\begin{displaymath}
\textrm{Let}\qquad y=10, \Big(4^2\cdot 10\Big)\alpha_1S_2=10\Rightarrow 9\alpha_1S_2=10.
\end{displaymath}
\begin{displaymath}
\textrm{Let}\qquad y=11, \Big(4^2\cdot 11\Big)\alpha_1S_2=11\Rightarrow 10\alpha_1S_2=11.
\end{displaymath}
Hence,
\begin{displaymath}
\alpha_1S_2=(0~1~2)(3~4~5)(6~7~8)(9~10~11)=\alpha^2=R_1.
\end{displaymath}
Consider
\begin{displaymath}
yx^2 =y\beta_1T_2
\end{displaymath}
and fix $x=4$.
\begin{displaymath}
\textrm{Let}\qquad y=0, 0\cdot 4^2=0\beta_1T_2\Rightarrow 2=0\beta_1T_2.
\end{displaymath}
\begin{displaymath}
\textrm{Let}\qquad y=1, 1\cdot 4^2=1\beta_1T_2\Rightarrow 0=1\beta_1T_2.
\end{displaymath}
\begin{displaymath}
\textrm{Let}\qquad y=2, 2\cdot 4^2=2\beta_1T_2\Rightarrow 1=2\beta_1T_2.
\end{displaymath}
\begin{displaymath}
\textrm{Let}\qquad y=3, 3\cdot 4^2=3\beta_1T_2\Rightarrow 5=3\beta_1T_2.
\end{displaymath}
\begin{displaymath}
\textrm{Let}\qquad y=4, 4\cdot 4^2=4\beta_1T_2\Rightarrow 3=4\beta_1T_2.
\end{displaymath}
\begin{displaymath}
\textrm{Let}\qquad y=5, 5\cdot 4^2=5\beta_1T_2\Rightarrow 4=5\beta_1T_2.
\end{displaymath}
\begin{displaymath}
\textrm{Let}\qquad y=6, 6\cdot 4^2=6\beta_1T_2\Rightarrow 8=6\beta_1T_2.
\end{displaymath}
\begin{displaymath}
\textrm{Let}\qquad y=7, 7\cdot 4^2=7\beta_1T_2\Rightarrow 6=7\beta_1T_2.
\end{displaymath}
\begin{displaymath}
\textrm{Let}\qquad y=8, 8\cdot 4^2=8\beta_1T_2\Rightarrow 7=8\beta_1T_2.
\end{displaymath}
\begin{displaymath}
\textrm{Let}\qquad y=9, 9\cdot 4^2=9\beta_1T_2\Rightarrow 11=9\beta_1T_2.
\end{displaymath}
\begin{displaymath}
\textrm{Let}\qquad y=10, 10\cdot 4^2=10\beta_1T_2\Rightarrow 9=10\beta_1T_2.
\end{displaymath}
\begin{displaymath}
\textrm{Let}\qquad y=11, 11\cdot 4^2=11\beta_1T_2\Rightarrow 10=11\beta_1T_2.
\end{displaymath}
Hence,
\begin{displaymath}
\beta_1T_2=(0~2~1)(3~5~4)(6~8~7)(9~11~10)=\alpha^{-2}=R_2.
\end{displaymath}
Consider
\begin{displaymath}
(x^2yx^{-2})\gamma_1\mathcal{R}_2 =y
\end{displaymath}
and fix $x=4$.
\begin{displaymath}
\textrm{Let}\qquad y=0, \Big(4^2\cdot 0\cdot
4^2\Big)\gamma_1\mathcal{R}_2=0\Rightarrow 0\gamma_1\mathcal{R}_2=0.
\end{displaymath}
\begin{displaymath}
\textrm{Let}\qquad y=1, \Big(4^2\cdot 1\cdot
4^2\Big)\gamma_1\mathcal{R}_2=1\Rightarrow 1\gamma_1\mathcal{R}_2=1.
\end{displaymath}
\begin{displaymath}
\textrm{Let}\qquad y=2, \Big(4^2\cdot 2\cdot
4^2\Big)\gamma_1\mathcal{R}_2=2\Rightarrow 2\gamma_1\mathcal{R}_2=2.
\end{displaymath}
\begin{displaymath}
\textrm{Let}\qquad y=3, \Big(4^2\cdot 3\cdot
4^2\Big)\gamma_1\mathcal{R}_2=3\Rightarrow 3\gamma_1\mathcal{R}_2=3.
\end{displaymath}
\begin{displaymath}
\textrm{Let}\qquad y=4, \Big(4^2\cdot 4\cdot
4^2\Big)\gamma_1\mathcal{R}_2=4\Rightarrow 4\gamma_1\mathcal{R}_2=4.
\end{displaymath}
\begin{displaymath}
\textrm{Let}\qquad y=5, \Big(4^2\cdot 5\cdot
4^2\Big)\gamma_1\mathcal{R}_2=5\Rightarrow 5\gamma_1\mathcal{R}_2=5.
\end{displaymath}
\begin{displaymath}
\textrm{Let}\qquad y=6, \Big(4^2\cdot 6\cdot
4^2\Big)\gamma_1\mathcal{R}_2=6\Rightarrow 6\gamma_1\mathcal{R}_2=6.
\end{displaymath}
\begin{displaymath}
\textrm{Let}\qquad y=7, \Big(4^2\cdot 7\cdot
4^2\Big)\gamma_1\mathcal{R}_2=7\Rightarrow 7\gamma_1\mathcal{R}_2=7.
\end{displaymath}
\begin{displaymath}
\textrm{Let}\qquad y=8, \Big(4^2\cdot 8\cdot
4^2\Big)\gamma_1\mathcal{R}_2=8\Rightarrow 8\gamma_1\mathcal{R}_2=8.
\end{displaymath}
\begin{displaymath}
\textrm{Let}\qquad y=9, \Big(4^2\cdot 9\cdot
4^2\Big)\gamma_1\mathcal{R}_2=9\Rightarrow 9\gamma_1\mathcal{R}_2=9.
\end{displaymath}
\begin{displaymath}
\textrm{Let}\qquad y=10, \Big(4^2\cdot 10\cdot
4^2\Big)\gamma_1\mathcal{R}_2=10\Rightarrow
10\gamma_1\mathcal{R}_2=10.
\end{displaymath}
\begin{displaymath}
\textrm{Let}\qquad y=11, \Big(4^2\cdot 11\cdot
4^2\Big)\gamma_1\mathcal{R}_2=11\Rightarrow
11\gamma_1\mathcal{R}_2=11.
\end{displaymath}
Hence,
\begin{displaymath}
\gamma_1\mathcal{R}_2=(0)(1)(2)(3)(4)(5)(6)(7)(8)(9)(10)(11)=I=R_0.
\end{displaymath}
So,
\begin{displaymath}
\alpha_1S_2=(0~1~2)(3~4~5)(6~7~8)(9~10~11)=\alpha^2=R_1,
\end{displaymath}
\begin{displaymath}
\beta_1T_2=(0~2~1)(3~5~4)(6~8~7)(9~11~10)=\alpha^{-2}=R_2,
\end{displaymath}
\begin{displaymath}
\gamma_1\mathcal{R}_2=(0)(1)(2)(3)(4)(5)(6)(7)(8)(9)(10)(11)=I=R_0.
\end{displaymath}
Therefore,
\begin{displaymath}
(\alpha_1S_2,\beta_1T_2,\gamma_1\mathcal{R}_2)=(\alpha^2,
\alpha^{-2}, I)=(R_1, R_2, R_0)=(R_1,R_2,I)\in Aut(L), \alpha
=R_{10}\in\Pi_\rho.
\end{displaymath}
This is a principal autotopism.
\paragraph{}
Also, we can construct an autotopism for the C-loop whose unbordered
multiplication table is in \cite{phi} by taking  the steps of the
construction above. Fix $x=9$ and let $\mu
=(0~13~5~14)(1~15~4~12)(2~9~10~8)(3~7~11~6)=R_{13}$. Then,
\begin{displaymath}
\alpha_1S_2=(0~5)(1~4)(2~10)(3~11)(6~7)(8~9)(12~15)(13~14)=\mu^2=R_5,
\end{displaymath}
\begin{displaymath}
\beta_1T_2=\alpha_1S_2=R_5,
\end{displaymath}
\begin{displaymath}
\gamma_1\mathcal{R}_2=I=R_0.
\end{displaymath}
Thus, $(\mu^2, \mu^2, I)=(R_5, R_5, R_0)\in Aut(L),
\mu=R_{13}\in\Pi_\rho $.

\end{document}